\documentclass[12pt,reqno,a4wide]{amsart}

\usepackage{tikz} 
\usetikzlibrary{positioning}
\usepackage{calrsfs}
\usepackage[charter]{mathdesign}

\allowdisplaybreaks

\title{ Non-decreasing Deutsch paths}

\author{Helmut Prodinger}
\address{Helmut Prodinger\\
	Department of Mathematical Sciences\\
	Stellenbosch University\\
	7602 Stellenbosch\\
	South Africa}
\email{hproding@sun.ac.za}

\begin{document}
\begin{abstract}
A variation of Dyck paths allows for down-steps of arbitrary length, not just one. This is motivated by ideas published by Emeric Deutsch around the turn of the millenium. We are interested in the subclass of them where the sequence of the levels of valleys is non-decreasing. This was studied around 20 years ago in the classical case. 
	\end{abstract}

\maketitle

\section{Introduction}

The paper \cite{Fezzi} introduced the subfamily of Dyck paths such that the level of the valleys is non-decreasing when scanning the path from left to right (``non-decreasing Dyck paths). The generating function
\begin{equation*}
\frac{1-2z}{1-3z+z^2}=\sum_{n\ge0}F_{2n-1}z^n,
\end{equation*}
where $2n$ denotes the length of the path, and $F_m$ are Fibonacci numbers, was already given. In \cite{elena}, the present writer found, that the usual translation of Dyck paths into plane trees, when restricted to the subfamily of non-decreasing trees, leads to a very simple tree structure, called Elena trees, from which it is easy to find generating functions, bijections, and consider various parameters of them. The paper \cite{DP} is also of relevance here. The study of various parameter was picked up again recently in \cite{Raimundo}.

The following two figures show such a path and the corresponding Elena tree.

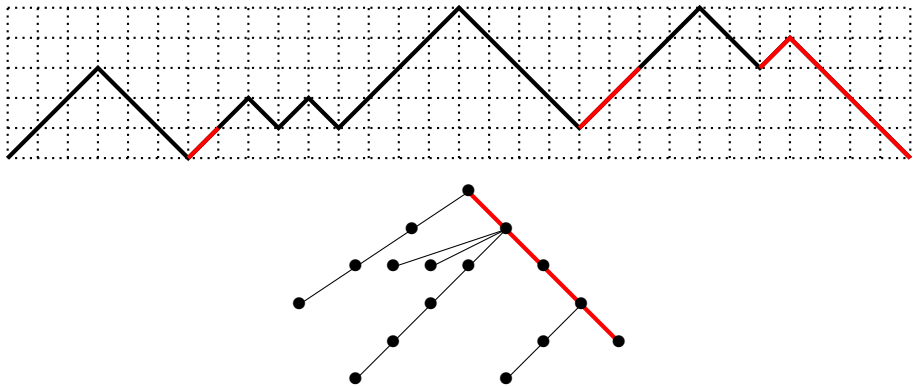
\begin{center}
		\begin{tikzpicture}[scale=0.4]
		
		\draw[step=1.cm,black,thick,dotted] (-0.0,-0.0) grid (30.0,5.0);

		\draw[ultra thick] (0,0) to (3,3) to (6,0) to (8,2) to (9,1) to (10,2) to (11,1) to (15,5) to (19,1) to (23,5)
		to (25,3)to (26,4)to (30,0);
		\draw[ultra thick,red] (6,0)to(7,1);
				\draw[ultra thick,red] (19,1)to(21,3);
				\draw[ultra thick,red] (25,3)to (26,4)to (30,0);
		\end{tikzpicture}
\end{center}
	
\begin{center}
	\begin{tikzpicture}[scale=0.5]
	
	\draw[ultra thick,red](0,0)to(4,-4);
	
	\node at (0,0) {$\bullet$};\node at (1,-1) {$\bullet$};\node at (2,-2) {$\bullet$};\node at (3,-3) {$\bullet$};
	\node at (4,-4) {$\bullet$};
	
	\draw(0,0)to(-4.5,-3);\node at (-1.5,-1) {$\bullet$};\node at (-3,-2) {$\bullet$};\node at (-4.5,-3) {$\bullet$};
	\draw(1,-1)to(-2,-2);\draw(1,-1)to(-1.0,-2);\draw(1,-1)to(-3,-5);
	\node at (-2,-2) {$\bullet$};\node at (-1,-2) {$\bullet$};\node at (0,-2) {$\bullet$};\node at (-1,-3) {$\bullet$};
	\node at (-2,-4) {$\bullet$};\node at (-3,-5) {$\bullet$};
	
	\draw(3,-3)to(1,-5);\node at (1,-5) {$\bullet$};\node at (2,-4) {$\bullet$};
	
	\end{tikzpicture}
	\end{center}
	When the Dyck path leaves a certain level  and never comes back to it (except for the final home-run), we draw a red edge.
	In the corresponding tree, this is the backbone, and the rest are just paths, hanging down.

	Recently, following ideas of Emeric Deutsch~\cite{Deutsch}, an extension of Dyck paths was studied in \cite{Deutsch-strip} (``Deutsch paths'').      
	Here, instead of just the usual down-step $(1,-1)$, all possible down-steps $(1,-j)$ for $j=1,2,3,\dots$ are allowed.
	
	The present note is dedicated to the description and characterization of non-decreasing Deutsch paths.

	\section{Enumeration }
	
	A sequence of up-steps followed by a sequence of down-steps ending at the same level is easy to enumerate for Dyck paths: 
	it exists only if the length is $2k$, and then there is just one such object.
	
	For Deutsch paths, this is a bit trickier, because of the various down-steps that are available and can be combined. We compute this number when the number of up-steps equals $j\ge1$ and the total number of steps is $k>j$:	
\begin{align*}
\sum_{1\le j< k}[z^{j}]\Big(\frac{z}{1-z}\Big)^{k-j}
&=\sum_{1\le j< k}[z^{2j-k}]\Big(\frac{1}{1-z}\Big)^{k-j}=\sum_{1\le j< k}\binom{j-1}{2j-k}
=\sum_{0\le j< k-1}\binom{k-j-2}{j}=F_{k-1}.
\end{align*}	
	The generating function of these numbers is
\begin{equation*}
\sum_{k\ge2}z^kF_{k-1}=\frac{z^2}{1-z-z^2}.
\end{equation*}
A sequence of such objects (or, equivalently, a bundle of such paths) has generating function
\begin{equation*}
\sum_{k\ge0}\Bigl(\frac{z^2}{1-z-z^2}\Bigr)^k=
\frac{1}{1-\dfrac{z^2}{1-z-z^2}}=\frac{1-z-z^2}{(1+z)(1-2z)}=:H.
\end{equation*}
This allows us to compute the total number of objects:
\begin{equation*}
\sum_{1\le j<k}H^j\binom{j-1}{k-1-j}z^k=
\sum_{1\le j}H^jz^{j+1}(1+z)^{j+1}=\frac {{z}^{2} ( 1-z-{z}^{2} ) }{ (1+z)(1-z) ( 1-2z-z^2 ) }.
\end{equation*}
We find it convenient to   allow the empty path as well, which means that we add 1 and get
\begin{equation*}
\frac{1}{4(1-z)}+\frac{1}{4(1+z)}+\frac12\frac{1-2z}{1-2z-z^2}=\cfrac{1}{1-\cfrac{z^2}{1-\cfrac{z}{1-\cfrac{z}{1-z^2}}}}.
\end{equation*}
We will explain this formula again in the next section in a more combinatorial fashion.
Now let $a=1+\sqrt2$, $b=1-\sqrt2$. Since
\begin{equation*}
	\frac12\frac{1-2z}{1-2z-z^2}=\frac12\frac{1-2z}{(1-az)(1-bz)}=
		\Bigl(\frac14-\frac18\sqrt2\Bigr)\frac1{1-az}+\Bigl(\frac14+\frac18\sqrt2\Bigr)\frac1{1-bz},
\end{equation*}
we found the number of non-decreasing Deutsch paths of length $n$:
\begin{equation*}
\frac14(1+(-1)^n)+\frac14(a^n+b^n)-\frac1{4\sqrt2}(a^n-b^n).
\end{equation*}
The numbers 
\begin{equation*}
\frac{a^n+b^n}2
\end{equation*}
are sequence A001333 in \cite{OEIS}, and the sequence  A000129 (Pell numbers)
\begin{equation*}
	\frac{a^n-b^n}{2\sqrt2}
\end{equation*}
is even more famous.

\section{Combinatorial considerations}

For increasing Dyck paths it is enough to count the number of up-steps, since the number of down-steps is the same.
This is no longer the case in the Deutsch model. So we have to count the total number of steps, which is the length of the (non-decreasing) Deutsch path.

When we form the associated tree, we must somehow indicate which downstep was used. One way of doing this is to use arrows to indicate where a new down-step starts.

To explain the concepts, we prepared a list of all non-decreasing Deutsch paths of length 5 and the corresponding trees.

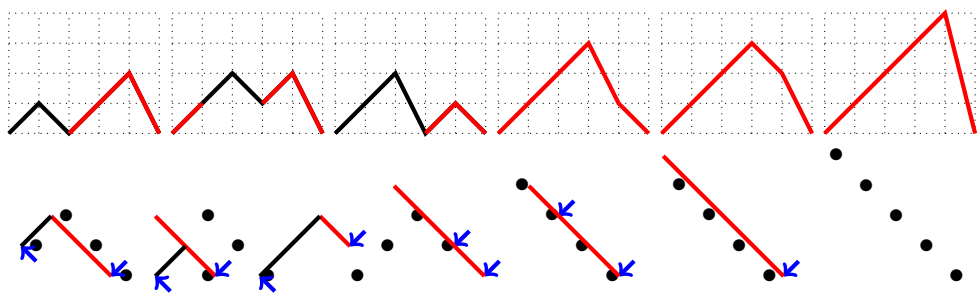
\begin{center}
\begin{tikzpicture}[scale=0.4]

\draw[step=1.cm,black,dotted] (-0.0,-0.0) grid (5.0,4.0);

\draw[ultra thick] (0,0) to (1,1) to (2,0) to (3,1) to (4,2) to (5,0);
\draw[ultra thick,red]  (2,0) to (3,1) to (4,2) to (5,0);
\end{tikzpicture}
\begin{tikzpicture}[scale=0.4]

\draw[step=1.cm,black,dotted] (-0.0,-0.0) grid (5.0,4.0);

\draw[ultra thick] (0,0) to (1,1) to (2,2) to (3,1) to (4,2) to (5,0);
\draw[ultra thick,red]  (0,0) to (1,1) ;
\draw[ultra thick,red]   (3,1) to (4,2) to (5,0);
\end{tikzpicture}
\begin{tikzpicture}[scale=0.4]

\draw[step=1.cm,black,dotted] (-0.0,-0.0) grid (5.0,4.0);

\draw[ultra thick] (0,0) to (1,1) to (2,2) to (3,0) to (4,1) to (5,0);
\draw[ultra thick,red]  (3,0) to (4,1) to (5,0);
\end{tikzpicture}
\begin{tikzpicture}[scale=0.4]

\draw[step=1.cm,black,dotted] (-0.0,-0.0) grid (5.0,4.0);

\draw[ultra thick,red] (0,0) to (1,1) to (2,2) to (3,3) to (4,1) to (5,0);

\end{tikzpicture}
\begin{tikzpicture}[scale=0.4]

\draw[step=1.cm,black,dotted] (-0.0,-0.0) grid (5.0,4.0);

\draw[ultra thick,red] (0,0) to (1,1) to (2,2) to (3,3) to (4,2) to (5,0);

\end{tikzpicture}
\begin{tikzpicture}[scale=0.4]

\draw[step=1.cm,black,dotted] (-0.0,-0.0) grid (5.0,4.0);

\draw[ultra thick,red] (0,0) to (1,1) to (2,2) to (3,3) to (4,4) to (5,0);

\end{tikzpicture}

\hspace*{-2.5cm}
\begin{tikzpicture}[scale=0.4]
\hspace*{0.2cm}
 \node at (0,0) {$\bullet$};
 \node at (1,-1) {$\bullet$};
  \node at (2,-2) {$\bullet$};
   \node at (-1,-1) {$\bullet$};

\draw[ultra thick,red] (0,0) to (2,-2);
\draw[ultra thick,blue,->] (2.5,-1.5) to (2,-2);
\draw[ultra thick,blue,->] (-0.5,-1.5) to (-1,-1);
\draw[ultra thick] (0,0) to (-1,-1);
\end{tikzpicture}
\begin{tikzpicture}[scale=0.4]
\hspace*{0.7cm}
\node at (0,0) {$\bullet$};
\node at (1,-1) {$\bullet$};
\node at (2,-2) {$\bullet$};
\node at (0,-2) {$\bullet$};

\draw[ultra thick,red] (0,0) to (2,-2);
\draw[ultra thick] (1,-1) to (0,-2);

\draw[ultra thick,blue,->] (2.5,-1.5) to (2,-2);
\draw[ultra thick,blue,->] (0.5,-2.5) to (0,-2);

\end{tikzpicture}
\begin{tikzpicture}[scale=0.4]
\hspace*{1.3cm}
\node at (0,0) {$\bullet$};
\node at (1,-1) {$\bullet$};
\node at (-2,-2) {$\bullet$};
\node at (-1,-1) {$\bullet$};

\draw[ultra thick,red] (0,0) to (1,-1);
\draw[ultra thick,blue,->] (1.5,-0.5) to (1,-1);
\draw[ultra thick,blue,->] (-1.5,-2.5) to (-2,-2);
\draw[ultra thick] (0,0) to (-2,-2);
\end{tikzpicture}
\begin{tikzpicture}[scale=0.4]
\hspace*{1.7cm}
\node at (0,0) {$\bullet$};
\node at (1,-1) {$\bullet$};
\node at (2,-2) {$\bullet$};
\node at (3,-3) {$\bullet$};

\draw[ultra thick,blue,->] (3.5,-2.5) to (3,-3);
\draw[ultra thick,blue,->] (2.5,-1.5) to (2,-2);
\draw[ultra thick,red] (0,0) to (3,-3);
\end{tikzpicture}
\begin{tikzpicture}[scale=0.4]
\hspace*{2.0cm}
\node at (0,0) {$\bullet$};
\node at (1,-1) {$\bullet$};
\node at (2,-2) {$\bullet$};
\node at (3,-3) {$\bullet$};

\draw[ultra thick,blue,->] (3.5,-2.5) to (3,-3);
\draw[ultra thick,blue,->] (1.5,-0.5) to (1,-1);
\draw[ultra thick,red] (0,0) to (3,-3);
\end{tikzpicture}
\begin{tikzpicture}[scale=0.4]
\hspace*{2.3cm}
\node at (0,0) {$\bullet$};
\node at (1,-1) {$\bullet$};
\node at (2,-2) {$\bullet$};
\node at (3,-3) {$\bullet$};
\node at (4,-4) {$\bullet$};

\draw[ultra thick,blue,->] (4.5,-3.5) to (4,-4);
\draw[ultra thick,red] (0,0) to (4,-4);
\end{tikzpicture}
\end{center}
We notice that the total number of edges plus the number of blue arrows is always 5, since each blue arrow indicates where a new down-step starts, and it must be counted just once. Alternatively, we can work with double edges, so these edges need to be counted twice. In the original Dyck model, every edge would be counted twice, leading to the length of the Dyck path.

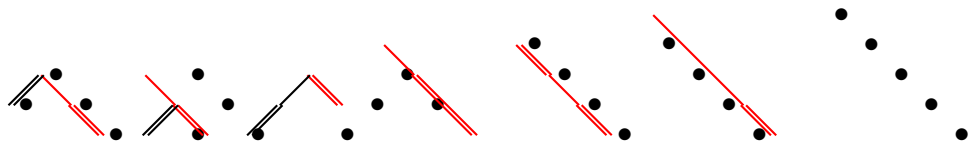
\begin{center}
	\hspace*{-2.5cm}
\begin{tikzpicture}[scale=0.4]
\hspace*{0.2cm}
\node at (0,0) {$\bullet$};
\node at (1,-1) {$\bullet$};
\node at (2,-2) {$\bullet$};
\node at (-1,-1) {$\bullet$};

\draw[ thick,red] (1,-1) to (0,0);

\draw[ thick,red,red,transform canvas={xshift=-1.0pt}] (1,-1) to (2,-2);
\draw[ thick,red,red,transform canvas={xshift=1.0pt}] (1,-1) to (2,-2);
\draw[ thick,transform canvas={xshift=-1.0pt}] (-1,-1) to (0,0);
\draw[ thick,transform canvas={xshift=1.0pt}] (-1,-1) to (0,0);

\end{tikzpicture}
\begin{tikzpicture}[scale=0.4]
\hspace*{0.7cm}
\node at (0,0) {$\bullet$};
\node at (1,-1) {$\bullet$};
\node at (2,-2) {$\bullet$};
\node at (0,-2) {$\bullet$};

\draw[ thick,red] (1,-1) to (0,0);

\draw[ thick,red,red,transform canvas={xshift=-1.0pt}] (1,-1) to (2,-2);
\draw[ thick,red,red,transform canvas={xshift=1.0pt}] (1,-1) to (2,-2);
\draw[ thick,transform canvas={xshift=-1.0pt}] (1,-1) to (0,-2);
\draw[ thick,transform canvas={xshift=1.0pt}] (1,-1) to (0,-2);


\end{tikzpicture}
\begin{tikzpicture}[scale=0.4]
\hspace*{1.3cm}
\node at (0,0) {$\bullet$};
\node at (1,-1) {$\bullet$};
\node at (-2,-2) {$\bullet$};
\node at (-1,-1) {$\bullet$};

\draw[  thick,red,transform canvas={xshift=-1.0pt}] (0,0) to (1,-1);
\draw[ thick,red,transform canvas={xshift=1.0pt}] (0,0) to (1,-1);
\draw[ thick,transform canvas={xshift=-1.0pt} ] (-1,-1) to (-2,-2);
\draw[ thick,transform canvas={xshift=1.0pt} ] (-1,-1) to (-2,-2);
\draw[ thick] (-1,-1) to (0,0);
\end{tikzpicture}
\begin{tikzpicture}[scale=0.4]
\hspace*{2.0cm}
\node at (0,0) {$\bullet$};
\node at (1,-1) {$\bullet$};
\node at (2,-2) {$\bullet$};
\node at (3,-3) {$\bullet$};

\draw[  thick,red ] (0,0) to (1,-1);
\draw[  thick,red,transform canvas={xshift=-1.0pt}] (3,-3) to (1,-1);
\draw[  thick,red,transform canvas={xshift=1.0pt}] (3,-3) to (1,-1);
\end{tikzpicture}
\begin{tikzpicture}[scale=0.4]
\hspace*{2.0cm}
\node at (0,0) {$\bullet$};
\node at (1,-1) {$\bullet$};
\node at (2,-2) {$\bullet$};
\node at (3,-3) {$\bullet$};

\draw[  thick,red,transform canvas={xshift=-1.0pt}] (0,0) to (1,-1);
\draw[  thick,red,transform canvas={xshift=1.0pt}] (0,0) to (1,-1);
\draw[  thick,red ] (2,-2) to (1,-1);
\draw[  thick,red,transform canvas={xshift=-1.0pt}] (3,-3) to (2,-2);
\draw[  thick,red,transform canvas={xshift=1.0pt}] (3,-3) to (2,-2);
\end{tikzpicture}
\begin{tikzpicture}[scale=0.4]
\hspace*{2.5cm}
\node at (0,0) {$\bullet$};
\node at (1,-1) {$\bullet$};
\node at (2,-2) {$\bullet$};
\node at (3,-3) {$\bullet$};
\node at (4,-4) {$\bullet$};

\draw[  thick,red] (3,-3) to (0,0);
\draw[  thick,red,transform canvas={xshift=-1.0pt}] (3,-3) to (4,-4);
\draw[  thick,red,transform canvas={xshift=1.0pt}] (3,-3) to (4,-4);

\end{tikzpicture}
\end{center}

We recall tilings with squares and
dominoes (length 2) from \cite{Q} and that length $n$ tilings are enumerated by the Fibonacci number $F_{n+1}$.

Note that $\frac{z}{1-z^2}$ enumerates sequences of type $122\dots2.$ Hence
\begin{equation*}
\frac1{1-\dfrac{z}{1-z^2}}=1+\sum_{n\ge1}[\text{tilings of length $n$, first tile is a square}]z^n.
\end{equation*}
Glueing one square to this starting square, the tilings of length $n+1$ starts with a domino.
If we delete the domino it is a general tiling of length $n-1$, enumerated by the Fibonacci number $F_n$.

So
\begin{equation*}
	\frac 1{1-\dfrac{z}{1-z^2}}=1+\sum_{n\ge1}F_nz^n\quad\text{and}\quad
	\frac z{1-\dfrac{z}{1-z^2}}=z+\sum_{n\ge2}F_{n-1}z^n.
\end{equation*}
That links the formula
\begin{equation*}
\frac{z^2}{1-z-z^2}=\sum_{n\ge2}F_{n-1}z^n
\end{equation*}
from before to  a tiling type argument. 

Using the description as a tree with single or double edges, we derive another form of the total generating function, which looks a little bit different:
\begin{equation*}
1+\cfrac{z^2}{1-\cfrac{z+z^2}{1-\cfrac{z^2}{1-z-z^2}}}\cdot\cfrac{1}{1-\cfrac{z^2}{1-z-z^2}}.
\end{equation*}

Now we explain the continued fraction 
\begin{equation*}
\cfrac{1}{1-\cfrac{z^2}{1-\cfrac{z}{1-\cfrac{z}{1-z^2}}}}
\end{equation*}
in a combinatorial way. We use tools that are common in combinatorics for words, in particular the symbolic equation
$(a+b)^*=b^*(ab^*)^*$, valid for any two nonempty formal languages or, power series that start with 0. We write $p$ for a path and $e$ for a single edge and $E$ for a double edge. As we discussed earlier, the generating function
\begin{equation*}
	\cfrac{z}{1-\cfrac{z}{1-z^2}}
\end{equation*}
corresponds to $e+p$. We split the backbone of the tree (the edges depicted in red) as $[(e\dots e)E]\dots[(e\dots e)E]$.
What hangs on the first node is enumerated by $b^*$.
Assume that the first edge of the backbone is $e$. Then the generating function according to this edge and the bundle of paths hanging down from node number 2 is enumerated by $ep^*$. Should the next edge also be $e$, it is enumerated via another $ep^*$, and so on. If the first group of single edges in the backbone has $l$ such edges, it is enumerated via $p^*(ep^*)^l$. The number $l$ can be any integer $\ge0$. However, the group ends with an $E$, and this one contributes $z^2$. This explains the term
\begin{equation*}
\Gamma=\cfrac{z^2}{1-\cfrac{z}{1-\cfrac{z}{1-z^2}}}.
\end{equation*}
But there is another `$*$' coming in, since the backbone consists of an arbitrary number of such groups. That explains the final formula
\begin{equation*}
\frac1{1-\Gamma}.
\end{equation*}

We switched freely between symbolic expressions (words/languages) and generating functions, where the variable $z$ counts an edge.
The book \cite{FS} describes very competently how this works.

\section{Conclusion}

We defined non-decreasing Deutsch paths and found the corresponding generating function. This was eased by translating it into a tree model using simple and double edges.  The explicit enumeration formula involves Pell numbers and (as an intermediate step) Fibonacci numbers. Parameters like the (average of the) number of double edges, degree of the root, length of the occurring paths etc.\ can be easily computed using a second variable in the generating function. The manipulations related to the rational functions that one gets in this way are best done by a computer. Since this is a routine procedure, we refrain from doing this here.

\end{document}